\documentclass[10pt,psamsfonts]{amsart}
\usepackage{amsmath}
\usepackage{amsthm}
\usepackage{amssymb}
\usepackage{amscd}
\usepackage{amsfonts}
\usepackage{amsbsy}
\usepackage{graphicx}

\newcommand{\R}{\ensuremath{\mathbb{R}}}

\newcommand{\p}{\partial}
\newcommand{\la}{\lambda}

\newcommand{\e}{\varepsilon}
\newcommand{\T}{\theta}

\newcommand{\rn}{{\mathbb{R}}^n}

\newcommand{\re}{\mbox{\normalfont{Re}}}

\newtheorem {theorem} {Theorem} 

\begin{document}

\title[Triple Hopf bifurcation in a tritrophic food chain model]
{Analytical study of a triple Hopf bifurcation in a tritrophic food chain model}

\author[J.-P. Fran\c coise and J. Llibre]
{Jean--Pierre Fran\c coise and Jaume Llibre}

\address{Universit\'e P.-M. Curie, Paris 6, Laboratoire Jacques--Louis Lions,
UMR 7598 CNRS, Site Chevaleret, 175 Rue du Chevaleret, 75013 Paris, France}
\email{Jean-Pierre.Francoise@upmc.fr}

\address{Departament de Matem\`{a}tiques, Universitat Aut\`{o}noma de Barcelona, 08193 Bellaterra, Barcelona, Catalonia, Spain}
\email{jllibre@mat.uab.cat}

\thanks{The first author is partially supported by an ANR grant ``Analyse non lin\'eaire et applications aux rythmes du vivant"
 $BLAN07-2-182920$. The second author is partially supported by a DGICYT
grant number MTM2005-06098-C02-01 and by a CICYT grant number
2005SGR 00550.}

\subjclass{34C37, 58F13, 92D25}

\keywords{limit cycle, Hopf bifurcation, population dynamics, averaging
theory}
\date{}
\dedicatory{}

\maketitle

\begin{abstract}
We provide an analytical proof of the existence of a stable periodic orbit contained in the region of coexistence of the three species of a tritrophic chain. The method used consists in analyzing a triple Hopf bifurcation. For some values of the parameters three limit cycles bear via this bifurcation. One is contained in the plane where the top--predator is absent. Another one is not contained in the domain of interest where all variables are positive. The third one is contained where the three species coexist. The techniques for proving these results have been introduced in previous articles (see \cite{BL1, Ll}) and are based on the averaging theory of second--order. Existence of this triple Hopf bifurcation has been previously discovered numerically in {\cite{KBK}}.
\end{abstract}

\section{Introduction}

During these last 80 years, after the seminal works of Lotka \cite{Lo}
and Volterra \cite{Vo}, one of the main topics in mathematical ecology has been the study of (di)trophic food chains. This has been made by analyzing
many different planar differential systems under the common name of prey--predator models, for instance see \cite{Ba, DDH}. The existence of limit cycles, attractors, and several kind of bifurcations are the characteristics of those models which have been used to explain the complex behaviors observed in such systems.

\smallskip

In the late seventies some interest in the mathematics of tritrophic food chain models (composed of prey, predator, and top--predator) appeared, see for example \cite{FW, Ga, FS} and Predator-Prey models with parasitic infection \cite{HF}. The model we analyze in this article describes a tritrophic food chain composed of a logistic prey (x), a Holling type II predator (y), and a Holling type II top--predator (z). After a rescaling of the variables, it is given by the following system of ordinary differential equations (see \cite{HP, MR2, KDR, KH, Mu} for more details):
\begin{equation}\label{L1}
\begin{array}{l}
\dot x= x \Big(\rho - \dfrac{x}{k} - \dfrac{a_1 y}{b_1 + x}\Big),\\
\\
\dot y= y \Big(\dfrac{a_1 x}{b_1 + x} - \dfrac{a_2 z}{b_2 + y}- d_1\Big),\\
\\
\dot z= z \Big(\dfrac{a_2 y}{b_2 + y}- d_2\Big).
\end{array}
\end{equation}
In order to preserve the biological meaning of the model, the $8$ parameters of this system are assumed to be strictly positive. Similar types of systems have been studied in the case when time scales of the variables are different so that methods of approximations of slow-fast systems can be applied \cite{Vi1, Vi2}. We emphasize that we do need here this type of approximation. To clarify possible applications of our result, we focus in a range where the population of the superpredator $z$ remains small (near $z=0$). So in our view here the superpredator remains of small amplitude compared with the two other populations. Within this range, we look for the possibility of a periodic rhythms which is stable and where the three species coexist.

\smallskip
The bifurcation analysis can be carried out with respect to the two parameters $(d_1,d_2)$. The normal form analysis around one of the stationary points (called latter $p_3$ in this article was first carried out by Klebanoff and Hastings (\cite{KH}) and latter improved by \cite{KR}. These authors (Kuznetsov and Rinaldi) also discovered numerically the existence of a strange attractor. In the limit of slow-fast systems Muratori and Rinaldi (\cite{MR2}) developed a singular perturbation approach which supports the existence of a homoclinic intersection of the stable and unstable manifolds of the limit cycle which is contained in the $(x,y)$-plane. Classical arguments in the theory of dynamical systems can then be used to deduce the existence of infinitely many limit cycles of increasing period (but nothing is known in general on their stability). See also (\cite{Vi1} and \cite{Vi2}).

\smallskip

We are interested in the limit cycles of system \eqref{L1}, mainly in the ones which come from a Hopf bifurcation. There are several papers dedicated to these limit cycles see for instance \cite{Ma, Ch, De, DH1, DH2, MR1, KDR}. But in all these papers the existence of a triple Hopf bifurcation was not proved analytically, the results there are essentially numerical. The existence of the triple Hopf bifurcation was discovered in {\cite{KBK}}. We prove here that there are systems \eqref{L1} having $1$ or $3$ small amplitude limit cycles coming from a Hopf bifurcation, and we show how to study the type of stability of such limit cycles. The tool for obtaining these results is the theory of averaging of second order. In fact as we shall show from the singular point bifurcate $3$ small amplitude limit cycles one of them is contained in $z>0$, the other in $z=0$ and the last in $z<0$. Of course the unique small amplitude limit cycles which have biological meaning are the ones contained in $z\geq 0$.
\smallskip

As a guide to the potential ``ecological user", we should emphasize that we can prove the existence of a Hopf bifurcation of second order provided that the logistic growth of the prey satisfies the two relations:
$$\rho = \dfrac{b_1 (a_1 + d_1)}{(a_1 - d_1) k},$$
and
$$k = \dfrac{2 a_1 b_1^2 d_1}{(a_1-d_1)^2 (a_1 b_1-2 b_2 d_1)},$$
and that the rate of exponential decay of the superpredator satifies:
$$d_2 = \dfrac{a_1 a_2 b_1^2}{ a_1^2 b_2 k + b_2 d_1^2 k + a_1 (b_1^2 - 2 b_2 d_1 k)}.$$
There is also a fourth condition which defines a full open set in the parameter space:
$$a_1b_1>2b_2d_1.$$

These four conditions ensure the existence of a periodic orbit contained in the domain of coexistence of the three species. A last inequality is required to ensure the stability of this periodic orbit. Altough from the viewpoint of applications imposing inequalities on the parameters is not constraining too much (provided they are compatible and we check that the set of existence of a stable periodic orbit is not empty by giving one example), it is of course much less natural to impose the three first conditions above. These conditions seem necessary for the analytical proof. Practically it might be enough to check that the values of the parameters remain close to those prescribed above but this last point is beyond the dynamical systems techniques we use here.

\smallskip

The paper is structured as follows. In the appendix  \ref{ap} we summarize the basic results on the averaging theory of second order that we shall need for studying the Hopf bifurcation. The explicit results on the Hopf bifurcation are stated in Theorem \ref{t1} at the end of section \ref{s2}. In this section is also proved Theorem \ref{t1}. An example showing the existence of three stable small amplitude limit cycles bifurcating from a singular point of system \eqref{L1} is given in Section 3.

\section{The Hopf bifurcation}\label{s2}

We separate the study of the Hopf bifurcation in different subsections.

\subsection{The singular point which will exhibit the Hopf bifurcation}

The differential system \eqref{L1} can have the following six singular points:
\[
\begin{array}{l}
p_1=(0,0,0),\\
\\
p_2=(k\, \rho,0,0),\\
\\
p_3=\left(\dfrac{b_1 d_1}{a_1-d_1}, -\dfrac{b_1 (b_1 d_1+(d_1-a_1) k\,  \rho)}{(a_1-d_1)^2 k}, 0\right),\\
\\
p_4=\left(0, \dfrac{b_2 d_2}{a_2-d_2}, -\dfrac{b_2 d_1}{a_2-d_2} \right),\\
\\
p_5=\left(\dfrac{A+\sqrt{B}}{2(a_2- d_2)}, \dfrac{b_2 d_2}{a_2-d_2}, \dfrac{b_2 (a_1-d_1) \sqrt{B}-b_2 C (a_2-d_2)}{(a_2-d_2)\left(\sqrt{B}+D\right)} \right),\\
\end{array}
\]
\[
\begin{array}{l}
p_6=\left(\dfrac{A-\sqrt{B}}{2(a_2- d_2)}, \dfrac{b_2 d_2}{a_2-d_2}, \dfrac{b_2 (a_1-d_1) \sqrt{B}+b_2 C (a_2-d_2)}{(a_2-d_2)\left(\sqrt{B}-D\right)} \right),
\end{array}
\]
where
\[
\begin{array}{l}
A= -a_2 b_1 + b_1 d_2 + a_2 k\,  \rho - d_2 k\,  \rho,\\
B=4 (a_2- d_2)(a_2(b_1+ k\,  \rho)^2 - d_2 (b_1^2 + 4 a_1 b_2 k + 2 b_1 k\,  \rho + k^2 \rho^2)),\\
C= a_1 b_1 + b_1 d_1 - a_1 k\,  \rho + d_1 k\,  \rho,\\
D= a_2 b_1 - b_1 d_2 + a_2 k\,  \rho - d_2 k\,  \rho.
\end{array}
\]
Of course such singular points exist always that their denominators are nonzero and the expression $B\ge 0$.

\smallskip

The study of the stability of $p_5$ and $p_6$ looks very tedious in full generality due to the very long expressions of their eigenvalues. So the analysis of their possible Hopf bifurcations is out of scope. Since the eigenvalues of $p_1$ and $p_2$ are always real, we must put our interest in the points $p_3$ and $p_4$. We choose here to study $p_3$ because it seems more meaningful to proceed in the neighborhood of $z=0$ if we keep in mind the situation of biological interest we presented in the introduction. It could be quite possible to proceed with $p_4$.

\smallskip

We shall study the Hopf bifurcation at the singular point $p_3$. The eigenvalues at this singular point are
\[
\begin{array}{ll}
\la_{\pm}= &\dfrac{1}{2 a_1 (a_1 - d_1) k} \left(-a_1 b_1 d_1 - b_1 d_1^2 + a_1 d_1 k\, \rho - d_1^2 k\, \rho \pm \sqrt{\Delta}\right),\\
& \\
\mu=& -d_2 + \dfrac{a_2 b_1 (b_1 d_1 + (-a_1 + d_1) k\, \rho)}{ b_1^2 d_1 - b_2 (a_1 - d_1)^2 k + b_1 (-a_1 + d_1) k\, \rho}
\end{array}
\]
where
\[
\Delta= d_1 (-4 a_1 (a_1 - d_1)^2 k (-b_1 d_1 + (a_1 - d_1) k\, \rho) +
   d_1 (a_1 (b_1 - k\, \rho) + d_1 (b_1 + k\, \rho))^2).
\]

It is well known that a necessary condition in order to have a Hopf bifurcation at $p_3$ is that their pair of complex eigenvalues $\la_{\pm}$ when $\Delta<0$ must cross the imaginary axis. So we will take
\begin{equation}\label{xx}
{\rm Re}(\la_{\pm})  = \dfrac{-a_1 b_1 d_1 - b_1 d_1^2 + a_1 d_1 k\, \rho - d_1^2 k\, \rho}{2 a_1 (a_1 - d_1) k}= \e^2 l,
\end{equation}
where $\e$ is a small parameter necessary for applying the averaging theory in order to study the Hopf bifurcation and $l$ is an arbitrary parameter. Initially we had taken instead of $\e^2 l$ the expression $\e l_1 +\e^2 l_2$. But applying the averaging theory of first order we do not obtain any information about the Hopf bifurcation and in order that the averaged function of first order $F_{10}$ (see the appendix) becomes identically zero we must take $l_1=0$. This is the reason that now we are taking directly $\e^2 l$ instead of $\e l_1 +\e^2 l_2$.

\smallskip

In order to apply the averaging theory we must write our differential system \eqref{L1} into the normal form of the averaging, i.e. into the form \eqref{42}. Then we need to take
\begin{equation}\label{yy}
\mu  = \e^2 m,
\end{equation}
where $m$ is an arbitrary parameter. As in the equation \eqref{xx} at the beginning we took $\e m_1 +\e^2 m_2$, but for the same reason than before we must take $m_1=0$. In short we are taking directly $\e^2 m$ instead of $\e m_1 +\e^2 m_2$.

\smallskip

Solving equalities \eqref{xx} and \eqref{yy} with respect to $d_2$ and $\rho$ we obtain that
\begin{equation}\label{zz}
\begin{array}{l}
d_2= \dfrac{a_1 a_2 b_1^2 d_1+ E \e^2+2 a_1 b_1 (d_1-a_1) k l
   m \e^4}{d_1 \left(b_2 k a_1^2+\left(b_1^2-2 b_2 d_1 k\right)
   a_1+b_2 d_1^2 k+ 2 a_1 b_1 (a_1-d_1) k l \e^2\right)},\\
\\
\rho= \dfrac{2 a_1 (a_1-d_1) k l \e^2+b_1 d_1 (a_1+d_1)}{(a_1-d_1)
   d_1 k},
\end{array}
\end{equation}
with $E= -b_2 k m d_1^3+a_1 \left(-m b_1^2-2 a_2 k l b_1+2 b_2 d_1 k m\right)
d_1+a_1^2 k (2 a_2 b_1 l-b_2 d_1 m)$. Now the eigenvalues of the singular point $p_3$ become
\[
\e^2 l\pm \dfrac{\sqrt{k \left(\e^2 k\, l \left(l \e^2-2 a_1+2
   d_1\right)-b_1 d_1\right)}}{k},\qquad \e^2 m.
\]

\subsection{The system in the normal form for applying the averaging theory}

One of the most difficult steps for applying the averaging theory to a given differential system, as our system \eqref{L1} in order to study its limit cycles, is to write it in the normal form \eqref{42} of the averaging theory. Moreover in our case we want to study the small amplitude limit cycles coming from a Hopf bifurcation. For doing all this we shall follow the next steps:

\medskip

\noindent(i) First we translate the singular point $p_3$ at the origin of coordinates doing the change of variables $(X,Y,Z)=(x,y,z)-p_3$.

\medskip

\noindent(ii) We write the linear part at the origin of the differential system $(\dot X, \dot Y, \dot Z)$ in its real Jordan normal form. For this we do the change of variables
\[
\left(
\begin{array}{c}
X\\
Y\\
Z
\end{array}
\right)=
\left(
\begin{array}{ccc}
-\dfrac{d_1 k}{F} & 0 & 1\\
 & & \\
-\dfrac{\e^2 k\, l}{F} & 1 & \dfrac{\e^2(2l-m)}{d_1}\\
 & & \\
0 & 0 & G
\end{array}
\right)
\left(
\begin{array}{c}
U\\
V\\
W
\end{array}
\right),
\]
where
\[
\begin{array}{l}
F= \sqrt{b_1 d_1-k \left(\e^2 k\, l \left(l \e^2-2 a_1+2
   d_1\right)\right)},\\
G=  H I/ \left[ a_1 a_2 b_1 d_1 k \left(2 (a_1-d_1) k l e^2+b_1 d_1\right) \right],\\
H= b_2 k d_1^3+a_1 (b_1^2-2e^2 k l b_1-2 b_2 d_1 k)
   d_1+a_1^2 k (2 b_1 l e^2+b_2 d_1),\\
I= k \left(m (m-2 l) e^2+2 a_1 l-2 d_1 l\right) e^2+b_1 d_1.
\end{array}
\]

\noindent(iii) We write the new differential system $(\dot U, \dot V, \dot W)$ in the cylindrical coordinates $(R, \T, W)$ defined through  $U= R \cos \T$, $V= R \sin \T$ and $W= W$.

\medskip

\noindent(iv) In order to study the small amplitude limit cycles around the origin of coordinates we do the rescaling $(R, \T, W)= (\e r, \T, \e w)$.

\medskip

\noindent(v) Finally the system $(\dot r, \dot \T, \dot w)$ is written as $(r',w')= (dr/d\T, dw/d\T)$ and we obtain the tritrophic food chain model in the normal form of averaging:
\[
\begin{array}{l}
r'= \e F_{11}(\T,r,w)+ \e^2 F_{21}(\T,r,w)+O(\e^3),\\
w'= \e F_{12}(\T,r,w)+ \e^2 F_{22}(\T,r,w)+O(\e^3),
\end{array}
\]
where
\[
\begin{array}{l}
F_{11}= R_1/T_0,\\
F_{21}= (R_2 T_0-R_1 T_1)/T_0^2,\\
F_{21}= W_1/T_0,\\
F_{22}= (W_2 T_0-W_1 T_1)/T_0^2,
\end{array}
\]
and
\[
\begin{array}{ll}
R_1=& \dfrac{d_1^2 r^2 }{a_1 \sqrt{b_1 d_1 k}}\cos ^3\T -\dfrac{(a_1-d_1) r^2}{b_1}\cos^2\T \sin \T -\dfrac{2 d_1 w r}{a_1 k}\cos^2\T +\\
& \\
& \dfrac{\sqrt{b_1 d_1 k} w^2}{a_1 k^2} \cos \T+
\dfrac{(a_1-d_1)}{a_1 \sqrt{b_1 d_1 k}\left(b_2 k a_1^2+b_1^2 a_1-2 b_2 d_1 k a_1+b_2 d_1^2 k\right)^2} (b_2^2 k^2 a_1^5+\\
& \\
& a_2 b_1 b_2 k^2 a_1^4-3 b_2^2 d_1 k^2 a_1^4+ 2 b_1^2 b_2 k a_1^4+b_1^4 a_1^3+2 b_2^2 d_1^2 k^2 a_1^3-3 a_2 b_1 b_2 d_1 k^2 a_1^3-\\
& \\
& 2 b_1^2 b_2 d_1 k a_1^3+2 b_2^2 d_1^3 k^2 a_1^2+ 3 a_2 b_1 b_2 d_1^2 k^2 a_1^2+b_1^4 d_1 a_1^2-2 b_1^2 b_2 d_1^2 k a_1^2-\\
& \\
& 3 b_2^2 d_1^4 k^2 a_1-a_2 b_1 b_2 d_1^3 k^2 a_1+2b_1^2 b_2 d_1^3 k a_1+b_2^2 d_1^5 k^2)r w  \cos \T \sin \T+\\
& \\
& \dfrac{(a_1-d_1) d_1 (d_1-a_1) k r^2}{a_1 b_1 \sqrt{b_1 d_1 k}} \cos \T \sin^2\T +\\
& \\
& \dfrac{b_1 (a_1-d_1)^2 r w }{b_2 k a_1^2+b_1^2 a_1-2 b_2 d_1 k a_1+b_2 d_1^2 k}\sin ^2\T- \dfrac{(a_1-d_1)w^2} {a_1 k}\sin \T,\\
& \\
& \\
R_2=& \dfrac{b_2 k^2 r^2 w (a_1-d_1)^6}{a_1 b_1 \left(b_2 k a_1^2+b_1^2 a_1-2 b_2 d_1 k a_1+b_2 d_1^2 k\right)^2}\sin ^3\T -\dfrac{d_1 r^3(a_1-d_1)^2}{a_1^2 b_1^2} \cos ^4\T +\\
& \\
& \dfrac{w^3 (a_1-d_1)^2}{a_1^2 b_1 k}\sin \T- \dfrac{r (w^2 a_1^3-3 d_1 w^2 a_1^2-b_1^2 l a_1^2+3 d_1^2 w^2 a_1-d_1^3 w^2)}{a_1^2 b_1^2}\sin ^2\T+\\
& \\
& \dfrac{(a_1-d_1) d_1 (d_1-a_1)k r^3}{a_1 b_1^2 \sqrt{b_1 d_1 k}}\cos ^3\T \sin \T+\dfrac{3 (a_1-d_1)^2 \sqrt{b_1 d_1 k} w r^2}{a_1^2 b_1^2 k}\cos^3\T-\\
& \\
& \dfrac{(a_1-d_1)^3 d_1 k r^3}{a_1^2 b_1^3}\cos^2\T \sin^2\T +
   \dfrac{(a_1-d_1)^2 (2a_1+d_1) w r^2}{a_1^2 b_1^2}\cos^2\T \sin \T+\\
& \\
&   \dfrac{(-3 w^2 a_1^2+b_1 k l a_1^2+6 d_1 w^2 a_1-3 d_1^2 w^2)r} {a_1^2 b_1 k}\cos ^2\T + \dfrac{J}{K}\cos \T \sin ^2\T+\\
& \\
& \dfrac{\sqrt{b_1 d_1 k} w^3 (a_1-d_1)^2}{a_1^2 b_1 d_1 k^2}\cos \T -
   \dfrac{(a_1+2 d_1) r w^2 (a_1-d_1)^2}{a_1^2 b_1 \sqrt{b_1
   d_1 k}}\cos \T \sin \T,\\
& \\
& \\
T_0=& -\dfrac{\sqrt{b_1 d_1 k}}{k},\\
&\\
& \\
\end{array}
\]
\[
\begin{array}{ll}
T_1=& -\dfrac{(a_1-d_1) d_1 r}{a_1  b_1} \cos ^3\T-
   \dfrac{2 (d_1-a_1) \sqrt{b_1 d_1 k} w}{a_1 b_1 k} \cos^2\T-\\
& \\
& \dfrac{\sqrt{b_1 d_1 k}(k a_1^2-2 d_1 k a_1+b_1 d_1+d_1^2 k) r}{a_1 b_1^2 k}\cos^2\T \sin \T- \dfrac{(a_1-d_1) w^2}{a_1 k r}\cos \T+\\
& \\
&  \dfrac{(b_1 k a_1^3-2 b_1 d_1 k a_1^2+2 b_2 d_1 k a_1^2+2 b_1^2
   d_1 a_1+b_1 d_1^2 k a_1-4 b_2 d_1^2 k a_1+2 b_2 d_1^3 k) w}
   {a_1 k (b_2 k a_1^2+b_1^2 a_1-2 b_2 d_1 k a_1+b_2 d_1^2 k)}\cos \T \sin \T+\\
& \\
& \dfrac{(a_1-d_1)^2 r} {a_1 b_1}\cos\T \sin^2\T - \dfrac{L}{M} \sin ^2\T-
   \dfrac{\sqrt{b_1 d_1 k} w^2}{a_1 k^2 r}\sin\T,\\
\end{array}
\]
\[
\begin{array}{ll}
& \\
W_1=& \dfrac{a_2 b_1 b_2 (a_1-d_1)^3(a_1+d_1) k w }{\left(b_2 k
   a_1^2+b_1^2 a_1-2 b_2 d_1 k a_1+b_2 d_1^2 k\right)^2}\sin \T,\\
& \\
& \\
W_2=& m w-\dfrac{a_2 b_1^2 b_2 (a_1-d_1)^4 (a_1+d_1)^2 k w} {(b_2 k
   a_1^2+b_1^2 a_1-2 b_2 d_1 k a_1+b_2 d_1^2 k)^3}\sin ^2\T.
\end{array}
\]
where
\[
\begin{array}{ll}
J=& k(2 b_2^3 k^3 d_1^7-12 a_1 b_2^3 k^3 d_1^6+30 a_1^2 b_2^3 k^3 d_1^5+6
   a_1 b_1^2 b_2^2 k^2 d_1^5-40 a_1^3 b_2^3 k^3 d_1^4-\\
& 24 a_1^2 b_1^2 b_2^2 k^2 d_1^4+30 a_1^4 b_2^3 k^3 d_1^3+36 a_1^3 b_1^2
   b_2^2 k^2 d_1^3+ a_1^2 a_2 b_1^3 b_2 k^2 d_1^3+ 6 a_1^2 b_1^4 b_2 k  d_1^3-\\
& 12 a_1^5 b_2^3 k^3 d_1^2-24 a_1^4 b_1^2 b_2^2 k^2 d_1^2-3 a_1^3 a_2 b_1^3 b_2 k^2 d_1^2-12 a_1^3 b_1^4 b_2 k d_1^2+2 a_1^3 b_1^6 d_1+\\
& 2 a_1^6 b_2^3 k^3 d_1+6 a_1^5 b_1^2 b_2^2 k^2 d_1+3 a_1^4 a_2 b_1^3 b_2 k^2 d_1+6 a_1^4 b_1^4 b_2 k d_1-a_1^5 a_2 b_1^3 b_2 k^2) r^2 w\\
& (a_1-d_1)^3,\\
& \\
K=& a_1^2 b_1^2 \sqrt{b_1 d_1 k}(b_2 k a_1^2+b_1^2 a_1-2 b_2 d_1 k
   a_1+b_2 d_1^2 k)^3,\\
& \\
L=& (a_1-d_1)^2 \sqrt{b_1 d_1 k} (b_2^2 k^2 a_1^4+a_2 b_1 b_2 k^2 a_1^3-4
   b_2^2 d_1 k^2 a_1^3+2 b_1^2 b_2 k a_1^3+b_1^4 a_1^2+\\
& 6 b_2^2 d_1^2 k^2 a_1^2-2 a_2 b_1 b_2 d_1 k^2 a_1^2-4 b_1^2 b_2 d_1 k a_1^2-4   b_2^2 d_1^3 k^2 a_1+a_2 b_1 b_2 d_1^2 k^2 a_1+\\
& 2 b_1^2 b_2 d_1^2 k a_1+b_2^2 d_1^4 k^2) w,\\
& \\
M=& a_1 b_1 d_1 k (b_2 k a_1^2+b_1^2 a_1-2 b_2 d_1 k a_1+b_2 d_1^2
   k)^2.
\end{array}
\]

\subsection{The computation of the small amplitude limit cycles}

Now we compute the function $F_{10}$, see Theorem \ref{ave2}, and we obtain
\[
\begin{array}{rl}
F_{10}(r,w)=& \displaystyle{\left(\dfrac{1}{2\pi}\int _0^{2\pi}F_{11}(\T,r,w)d\T, \dfrac{1}{2\pi}\int _0^{2\pi}F_{21}(\T,r,w)d\T \right)}\\
& \\
=& \left(0, -\dfrac{N w}{P}\right),
\end{array}
\]
where
\[
\begin{array}{ll}
N=& b_1 k a_1^3-2 b_1 d_1 k a_1^2-2 b_2 d_1 k a_1^2-2 b_1^2
   d_1 a_1+b_1 d_1^2 k a_1+4 b_2 d_1^2 k a_1-2 b_2 d_1^3 k,\\
& \\
P=& a_1 b_2 k a_1^2+b_1^2 a_1-2 b_2 d_1 k a_1+b_2 d_1^2 k).
\end{array}
\]
Taking
\begin{equation}\label{L2}
k= \dfrac{2 a_1 b_1^2 d_1}{(a_1-d_1)^2 (a_1 b_1-2 b_2 d_1)},
\end{equation}
we obtain that the averaged function of first order is identically zero. We must compute the averaged function of second order, for more details see the appendix.

\smallskip

We note that since in the expressions of $F_{ij}$ with $i,j\in\{1,2\}$ appears $\sqrt{k}$ we need that
\begin{equation}\label{L3}
a_1 b_1-2 b_2 d_1> 0.
\end{equation}

{From} \eqref{f2} the averaged function of second order $F_{20}(r,w)$ has the two components:
\[
\begin{array}{ll}
F_{201}(r,w)=& \dfrac{N r}{2 \sqrt{2} a_1^4 b_1^4 d_1}
\big(2(a_1-d_1)^2 (a_1 b_1-2 b_2 d_1)(b_1 a_1^3-b_1 d_1 a_1^2-4 a_2 b_2 d_1^2) w^2-\\
& \\
& a_1^3 b_1^2 d_1 (4 a_1 l b_1^2+ d_1(d_1 -a_1) r^2)\big),\\
& \\
\end{array}
\]
\[
\begin{array}{rl}
F_{202}(r,w)=& -\dfrac{\sqrt{2}}{a_1^5 b_1^5} (a_1-d_1)^2 (a_1 b_1-2
   b_2 d_1) N^3 w \\
& \\
& \left(a_1^4 m_2 b_1^4-6 a_1 a_2 b_2 d_1^3 r^2 b_1-2 a_2 b_2 (a_1-d_1)^2
   d_1 (a_1 b_1-2 b_2 d_1)w^2\right),
\end{array}
\]
where
\[
N= \sqrt{\dfrac{a_1 b_1}{(a_1-d_1)^2 (a_1 b_1-2 b_2 d_1)}}.
\]

In order to look for the small amplitude limit cycles bifurcating from the origin of system \eqref{L1}, after all the changes of coordinates that we did and according with Theorem \ref{ave2}, we must find the zeros $(r_0,w_0)$ with $r_0>0$ of the system
\begin{equation}\label{L4}
F_{201}(r,w)=0, \qquad F_{202}(r,w)=0.
\end{equation}
such that the Jacobian
\begin{equation}\label{L5}
\left. \det\left(
\begin{array}{cc}
\dfrac{\p F_{201}}{\p r} & \dfrac{\p F_{201}}{\p w}\\
\\
\dfrac{\p F_{202}}{\p r} & \dfrac{\p F_{202}}{\p w}
\end{array}
\right)\right|_{(r,w)= (r_0,w_0)}
\end{equation}
be nonzero.

\smallskip

It is easy to check that system \eqref{L4} can have at most $3$ solutions satisfying \eqref{L5} according with the values of the parameters of system \eqref{L1}. More precisely one solution $(r_1,w_1)$ is
\begin{equation}\label{z1}
r_1= 2b_1 \sqrt{\dfrac{a_1 l}{(a_1-d_1)d_1}}= 2b_1 \sqrt{R_1}, \qquad w_1=0;
\end{equation}
of course it exists if $R_1>0$. The other two possible solutions are $(r_2, \pm w_2)$ where
\begin{eqnarray}
r_2 &=& \dfrac{a_1 b_1}{d_1} \sqrt{\dfrac{a_1 b_1 \left(a_1^2 b_1 (a_1-d_1) m-4 a_2 b_2 d_1^2 (l+m)\right)}{a_2 b_2 d_1 \left(5 b_1 a_1^3-5 b_1 d_1 a_1^2-24 a_2 b_2 d_1^2\right)}}= \dfrac{a_1 b_1}{d_1} \sqrt{R_2},\nonumber\\
& &  \label{z2}\\
w_2& =& \dfrac{a_1^2 b_1^2}{\sqrt{2}} \sqrt{\dfrac{24 a_2 b_2 d_1^2 l-a_1^2 b_1 (a_1-d_1)m}{a_2 b_2 (a_1-d_1)^2 d_1 (2b_2 d_1-a_1 b_1) \left(-5 b_1 a_1^3+5 b_1 d_1 a_1^2+24 a_2 b_2 d_1^2\right)}}\nonumber\\
&=& \dfrac{a_1^2 b_1^2}{\sqrt{2}} \sqrt{W_2};\nonumber
\end{eqnarray}
of course again these last two small amplitude limit cycles will exist if $R_2>0$ and $W_2>0$.

\smallskip

When the parameters of system \eqref{L1} are such that we have $3$ small amplitude limit cycles, then the one which has initial conditions on the plane $w=0$ (see Theorem \ref{ave2}) remains in this plane because it is invariant by the flow of the system. Again since this plane is invariant, the other two small amplitude limit cycles which have initial conditions on half--spaces $w>0$ and $w>0$ remain in such half--spaces. Of course the invariant plane $w=0$ corresponds to the plane $z=0$ in the initial coordinates.

\smallskip

According with the statement (b) of Theorem \ref{ave2} we can compute the type of stability of these small amplitude limit cycles computing the eigenvalues of the matrix \eqref{L5}. More precisely, let $\la_1$ and $\la_2$ be the two eigenvalues of the matrix \eqref{L5} evaluated on a zero $(r_0,w_0)$ of system \eqref{L4}, then the small amplitude limit cycle associated to the zero $(r_0,w_0)$
\begin{itemize}
\item[(I)] is a local repeller if $\re(\la_1), \re(\la_2)>0$,
\item[(II)] is a local attractor if $\re(\la_1), \re(\la_2)<0$,
\item[(III)] has two invariant manifolds, one stable and the other unstable, which locally are formed by two $2$--dimensional cylinders.
\end{itemize}

We note that conditions \eqref{zz} with $\e=0$, \eqref{L2} and \eqref{L3} are necessary in order to apply the second order averaging theory. More precisely we need that system \eqref{L1} satisfies the conditions
\begin{equation}\label{zt}
\begin{array}{ll}
d_2 =& \dfrac{a_1 a2 b_1^2 d_1}{ a_1^2 b_2 d_1 k + b_2 d_1^3 k + a_1 d_1 (b_1^2 - 2 b_2 d_1 k)},\\
& \\
\rho =& \dfrac{b_1 (a_1 + d_1)}{(a_1 - d_1) k},\\
& \\
k =& \dfrac{2 a_1 b_1^2 d_1}{(a_1-d_1)^2 (a_1 b_1-2 b_2 d_1)},\\
& \\
0 <& a_1 b_1-2 b_2 d_1,
\end{array}
\end{equation}
in order that we can apply the averaging theory of second order for studying its small amplitude limit cycles.

\smallskip

In short we have proved the next result.

\begin{theorem}\label{t1}
The following statements hold.
\begin{itemize}
\item[(a)] The tritrophic food chain model given by system \eqref{L1} satisfying the four conditions \eqref{zt} has a triple Hopf bifurcation at the singular point  $p_3$ if $R_1$, $R_2$ and $W_2$ are positive (these last three expressions are defined in \eqref{z1} and \eqref{z2}).

\item[(b)] More precisely, under the assumptions of statement (a) three small amplitude limit cycles bifurcate from $p_3$ with initial conditions in the coordinates $(U,V,W)$ given by $(\e r_1,0,\e w_1)$ and $(\e r_2,0,\pm \e w_2)$, where the values of $r_i$ and $w_i$ are given in \eqref{z1} and \eqref{z2}. Moreover in the variables $(x,y,z)$ the small amplitude limit cycle coming from the initial conditions $(\e r_1,0,\e w_1)$ lies on the plane $z=0$, the one coming from the initial conditions $(\e r_2,0,\e w_2)$ lies in the half--space $z>0$, and the remainder one lies in the half--space $z<0$.

\item[(c)] Under the assumptions of statement (a) the kind of stability of the small amplitude limit cycles is determined by the eigenvalues of the matrix \eqref{L5} as it is explained in (I), (II) and (III).
\end{itemize}
\end{theorem}

\section{An example}

When $\e=0$ the relations \eqref{zt} for the values
\[
a_1 = 5, \quad a_2 = 0.1, \quad  b_1 = 3, \quad b_2 = 2, \quad d_1 = 0.4,
\]
become
\[
d_2= 0.09,\quad \rho= 27.74,\quad k= 0.13,\quad a_1 b_1-2 b_2 d_1= 13.4.
\]
These values are compatible with the biological conditions. Moreover the parameters $l$ and $m$ are free. Taking $l=400$ and $m=1$ we obtain that the eigenvalues of the matrix in the expression \eqref{L5} are
\[
\begin{array}{lllll}
-154.96  & \mbox{and} & -0.32   &  \mbox{for} &  (r_1,w_1)=(221.16,0),\\
-135.14  & \mbox{and} & -0.29   &  \mbox{for} &  (r_2,w_2)=(207.24,39),\\
-135.14  & \mbox{and} & -0.29   &  \mbox{for} &  (r_2,w_2)=(207.24,-39).
\end{array}
\]
Therefore the three small amplitude limit cycles are local attractors for these values of the parameters.

\section{The appendix: Averaging theory of second order}\label{ap}

We shall use the following result.

\begin{theorem}[Second order averaging method]\label{ave2}
We consider the following differential system
\begin{equation} \label{42}
\frac{dx}{dt} =\e F_{1}(t,x)+\e^2 F_{2}(t,x)+\e^3 R(t,x,\e),
\end{equation}
where $F_{1}, F_{2}:\mathbb{R}\times D\to \rn$, $R:\mathbb{R}\times
D\times (-\e_f,\e _f)\to \rn$ are continuous functions,
$T$--periodic in the first variable, and $D$ is an open subset of
$\rn$. We assume that
\begin{itemize}
\item[(i)] $F_{1}(t,\cdot)\in C^1(D)$ for all $t\in \mathbb{R}$, $F_{1}$, $F_{2}$, $R$ and $D_{x}F_{1}$ are locally Lipschitz with respect to
$x$, and $R$ is differentiable with respect to $\e$.\\
We define $F_{10}$, $F_{20}:D\rightarrow \mathbb{R}^{n}$ as

\begin{equation*}
F_{10}(z)= \dfrac{1}{T}\int _0^{T}F_{1}(s, z)ds,
\end{equation*}

\begin{equation}\label{f2}
F_{20}(z)=  \dfrac{1}{T}\int _0^{T}\left[D_z F_{1}(s, z) \cdot \int
_0^{s}F_1(t,z)dt + F_2(s,z)\right]ds,
\end{equation}

and assume moreover that
\item[(ii)] for $V\subset D$ an open and bounded set and for each $\e \in (-\e_f,\e_f)\setminus \{0\}$, there exists $a_\e \in V$ such that  $F_{10}(a_{\e})+\e F_{20}(a_{\e})=0$
and $d_B(F_{10}+\e F_{20},V,a_\e)\neq 0$.
\end{itemize}
Then the following statements hold.
\begin{itemize}
\item[(a)] For $|\e|>0 $ sufficiently small there exists a $T$--periodic solution $\varphi (t,\e)$ of system \eqref{42} such that $\varphi (0, \e)= a_\e$.

\item[(b)] If the function $F_{10}+\e F_{20}$ is $C^1$ the stability or instability of the limit cycle $\varphi (t, \e)$ is given by the stability or instability of the singular point $a_\e$ of the averaged system
\[
\frac{dx}{dt} =\e F_{10}(x)+\e^2 F_{20}(x),
\]
corresponding to system \eqref{42}. In fact, the singular point $a_\e$ of the averaged system has the stability behavior of the Poincar\'{e} map associated to the limit cycle $\phi(t,\e)$.
\end{itemize}
\end{theorem}

The proof of statement (a) can be found in \cite{BL1}, and of statement (b) is standard, see for instance \cite{V, SV, Ll}.

\smallskip

The expression $d_B(F_{10}+\e F_{20},V,a_\e)\neq 0$ means that the Brouwer degree of the function $F_{10}+\e F_{20}: V\to \R^n$ at the fixed point $a_\e$ is not zero. For instance it is sufficient that the Jacobian of the function $F_{10}+\e F_{20}$ at $a_\e$ is not zero in order that the mentioned Brouwer degree be nonzero.

\smallskip

If $F_{10}$ is not identically zero, then the zeros of $F_{10}+\e F_{20}$ are mainly the zeros of $F_{10}$ for $\e$ sufficiently small. In this case the previous result provides the {\it averaging theory of first order}.

\smallskip

If $F_{10}$ is identically zero and $F_{20}$ is not identically zero, then the zeros of $F_{10}+\e F_{20}$ are mainly the zeros of $F_{20}$ for $\e$ sufficiently small. In this case the previous result provides the {\it averaging theory of second order}.

\end{document}